\documentclass[12pt,english,reqno]{amsart}  

\usepackage[T1]{fontenc} 
\usepackage{geometry} 
\usepackage{amssymb,amsmath,amsthm} 

\usepackage{pdfpages}
\usepackage{graphicx}
\usepackage{tikz} 
\usetikzlibrary{arrows,automata} 

\usepackage{algorithm}
\usepackage{algpseudocode}

\usepackage{float}

\usepackage{hyperref}
\hypersetup{
    colorlinks=true,
    linkcolor=blue,
    filecolor=magenta,      
    urlcolor=cyan,
}

\usepackage{float} 

\newtheorem{theorem}{Theorem}[section]

\newtheorem{lemma}[theorem]{Lemma}
\newtheorem{proposition}[theorem]{Proposition}

\DeclareMathOperator{\Err}{\mathrm{Err}}
\DeclareMathOperator{\best}{\mathrm{best}}
\DeclareMathOperator{\last}{\mathrm{last}}

\title[Gradient-Free Stochastic Method for the 4- and 6-Body Problems]
{Computing Families of Periodic Solutions of the 4- and 6-Body Problems
Using an Adaptive Gradient-Free Stochastic Method}
\keywords{n-body problem, periodic solutions,  gradient-free methods, continuation methods}
\author{Oscar Perdomo}
\address{Department of Mathematical Sciences, Central Connecticut State University, New Britain, CT, USA}
\email{perdomoosm@ccsu.edu}
\date{\today}

\graphicspath{{images-ver2/}}

\begin{document}

\begin{abstract}
In this paper, we numerically compute planar periodic motions of the
4- and 6-body problems. We consider a 4-body configuration formed by two
opposite pairs of bodies and a 6-body configuration formed by two equilateral
triangles. More generally, for two regular $k$-gons with masses $1$
and $m_2$ in the respective groups, we derive a system of four second-order
differential equations for two radial and two angular variables; for $k=2$,
the groups are opposite pairs. With $r_1(0)=1$, the periodicity conditions lead to an
overdetermined nonlinear system consisting of eight return equations in five
unknowns. We show that, for collision-free solutions, four symmetry
conditions at $H=T/2$ imply these eight equations. For prescribed mass ratio and angular advance,
we solve the resulting four equations in four unknowns using Newton's
method. An adaptive stochastic black-box procedure, which uses only
function evaluations, provides initial approximations. The method searches
for solutions by sampling in an anisotropic box whose dimensions adapt using
information from previous successful iterations. Exact satisfaction of the
return conditions gives a motion that repeats up to rotation and relabeling;
an angular advance that is a rational multiple of $2\pi$ then gives a
periodic motion. We report numerical examples and check the return
residuals of the printed initial conditions using both the reduced equations
and independent integrations of the full Cartesian equations.
\end{abstract}

\maketitle

\section{Introduction}\label{intro}

In this paper, we numerically compute families of planar periodic solutions
of the $n$-body problem for $n=4$ and $n=6$. In both cases, we use the
symmetries of the configurations to reduce the equations of motion to a
closed system of four second-order differential equations for two radial
variables and two angular variables. We then impose return conditions on
the solutions of these reduced equations and use an adaptive gradient-free
stochastic method to find initial approximations to the resulting nonlinear
systems. We then use Newton's method to compute the initial conditions
reported here.

In order to describe both configurations using the same notation, let
$k\geq 2$ be an integer and consider two groups of $k$ bodies. We assume
that the bodies in the first group have mass $m_1>0$ and those in the second
group have mass $m_2>0$. Let $r_1(t)$ and $r_2(t)$ be positive functions
giving the distances of the bodies in each group from the origin, and let
$\theta(t)$ and $\beta(t)$ be angular variables. For convenience, let us
define
\[
e(\alpha)=(\cos\alpha,\sin\alpha).
\]
We assume that the bodies move according to the following ansatz:
\begin{eqnarray}\label{generalansatz}
{\bf r}_{j+1}(t)&=&r_1(t)e\left(\theta(t)+\frac{2\pi j}{k}\right),\nonumber\\
{\bf r}_{j+k+1}(t)&=&r_2(t)e\left(\beta(t)+\frac{\pi}{k}
+\frac{2\pi j}{k}\right),\qquad j=0,\ldots,k-1.
\end{eqnarray}
Thus, for $k\geq 3$, each group forms a regular $k$-gon with center at
the origin. The two groups are allowed to change their relative size and
orientation during the motion. Notice that the polar angle of the first
body in the second group is $\beta(t)+\pi/k$.

Without loss of generality, we can assume that $r_1(0)=1$ by changing the
unit of distance and that $m_1=1$ by changing the unit of mass. After fixing
these two units, we can choose the unit of time so that the gravitational
constant is $G=1$. From now on, $m_2$ denotes the mass of each body in the
second group in these units.

We consider the initial conditions
\begin{equation}\label{ic}
\begin{gathered}
r_1(0)=1,\qquad r_2(0)=x_2>0,\qquad
\theta(0)=\beta(0)=0,\\
\dot r_1(0)=\dot r_2(0)=0,\qquad
\dot\theta(0)=x_3,\qquad \dot\beta(0)=x_4.
\end{gathered}
\end{equation}
In particular, the two groups initially have angular separation $\pi/k$.

Thus, although $2k$ bodies are moving in the plane, the symmetry of the
ansatz allows us to describe the motion using only the four functions
$r_1$, $r_2$, $\theta$, and $\beta$. The corresponding reduced equations
of motion are derived in Section~\ref{DE} for every integer $k\geq2$,
with mass $1$ in the first group and mass $m_2>0$ in the second.
We then write the cases $k=2$ and $k=3$ explicitly.

For a prescribed angular advance $\theta_1$, we look for initial data
$x_2>0$, $x_3$, $x_4$, a mass $m_2>0$, and a time $T>0$ satisfying
\begin{equation}\label{systemneqgeneral}
\begin{gathered}
r_1(T)=1,\qquad r_2(T)=x_2,\qquad
\dot r_1(T)=\dot r_2(T)=0,\\
\dot\theta(T)=x_3,\qquad \dot\beta(T)=x_4,\\
\theta(T)-\beta(T)=\frac{2\pi}{k},\qquad \theta(T)=\theta_1.
\end{gathered}
\end{equation}
These are eight equations in the five unknowns $x_2,x_3,x_4,m_2,T$.
They need not be independent; the midpoint formulation below gives a
square system of sufficient conditions.
Assuming that the solution has no collisions on $[0,T]$, let us explain
why these equations are sufficient for a pseudo-periodic motion.
Here we call a motion pseudo-periodic if, after a fixed time, its
positions and velocities repeat up to a fixed rotation and a relabeling
of bodies with the same mass.

Let $R_\varphi$ be rotation about the origin through the angle $\varphi$,
that is,
\[
R_\varphi=
\begin{pmatrix}
\cos\varphi&-\sin\varphi\\
\sin\varphi&\cos\varphi
\end{pmatrix}.
\]
For convenience, write ${\bf a}_j={\bf r}_{j+1}$ and
${\bf b}_j={\bf r}_{j+k+1}$, taking the indices $j$ modulo $k$.
In particular, ${\bf b}_{-1}={\bf b}_{k-1}$. The angular conditions
in~\eqref{systemneqgeneral} give
\[
\theta(T)=\theta_1,\qquad
\beta(T)=\theta_1-\frac{2\pi}{k}.
\]
Using the ansatz~\eqref{generalansatz} and the return of the radii,
we obtain
\begin{equation}\label{positionreturngeneral}
\begin{aligned}
{\bf a}_j(T)
 &=e\left(\theta_1+\frac{2\pi j}{k}\right)
   =R_{\theta_1}{\bf a}_j(0),\\
{\bf b}_j(T)
 &=x_2e\left(\theta_1+\frac{\pi}{k}
                       +\frac{2\pi(j-1)}{k}\right)
   =R_{\theta_1}{\bf b}_{j-1}(0).
\end{aligned}
\end{equation}
Thus, rotating the configuration at $T$ by $-\theta_1$ brings the
first group back to its initial positions. The second group also
occupies its initial positions, with a cyclic change of labels.

We must also check the velocities. Let $J(x,y)=(-y,x)$, so that
$Je(\alpha)$ is the unit tangent vector to the circle at $e(\alpha)$.
Since the radial velocities vanish at both $0$ and $T$, and the angular
velocities return to $x_3$ and $x_4$, respectively, we have
\[
\begin{aligned}
\dot{\bf a}_j(0)&=x_3J{\bf a}_j(0),&
\dot{\bf a}_j(T)&=x_3J{\bf a}_j(T),\\
\dot{\bf b}_j(0)&=x_4J{\bf b}_j(0),&
\dot{\bf b}_j(T)&=x_4J{\bf b}_j(T).
\end{aligned}
\]
Rotations commute with $J$. It follows from
\eqref{positionreturngeneral} that
\begin{equation}\label{velocityreturngeneral}
\dot{\bf a}_j(T)=R_{\theta_1}\dot{\bf a}_j(0),\qquad
\dot{\bf b}_j(T)=R_{\theta_1}\dot{\bf b}_{j-1}(0).
\end{equation}
Therefore both positions and velocities agree with the initial ones
after the same rotation and relabeling.

Rotations preserve
distances, and a relabeling within either group preserves the masses.
Consequently these transformations preserve Newton's equations.
Uniqueness for the initial value problem then allows us to continue
the motion by rotated and relabeled copies of the first arc, with
\begin{equation}\label{pseudoperiodicreturn}
{\bf a}_j(t+T)=R_{\theta_1}{\bf a}_j(t),\qquad
{\bf b}_j(t+T)=R_{\theta_1}{\bf b}_{j-1}(t).
\end{equation}
This is the pseudo-periodicity that we seek. If $\theta_1$ is a rational
multiple of $2\pi$, we can choose a positive integer $M$, divisible by
$k$, such that $M\theta_1$ is a multiple of $2\pi$. After $M$ returns,
both the rotation and the relabeling are the identity. The labeled
motion is then periodic with period $MT$, which need not be its smallest
period.

Figure~\ref{fig:pseudoperiodicreturn} illustrates this argument for
$k=2$ and $k=3$. For $k=2$, rotating back exchanges bodies $3$ and $4$.
For $k=3$, the rotated bodies $4$, $5$, and $6$ occupy the initial
states of bodies $6$, $4$, and $5$, respectively.

\begin{figure}[htbp]
\centering
\resizebox{\linewidth}{!}{
\begingroup
\definecolor{ppred}{RGB}{194,48,54}
\definecolor{ppink}{RGB}{40,48,60}
\definecolor{ppblue}{RGB}{47,91,140}
\begin{tikzpicture}[x=1cm,y=1cm,>=stealth,
  every node/.style={font=\small},line cap=round,line join=round]
\foreach \ppk/\ppy in {2/0,3/-5.70}{%
  \begin{scope}[shift={(0,\ppy)}]
    \ifnum\ppk=2
      \node[anchor=west,font=\bfseries\small] at (-1.65,2.65)
        {$k=2$: two opposite pairs};
    \else
      \node[anchor=west,font=\bfseries\small] at (-1.65,2.65)
        {$k=3$: two equilateral triangles};
    \fi
    \foreach \ppx/\ppphase/\ppoffset/\pptitle in
      {0/0/1/{Initial state},4.7/38/-1/{State at $T$},
       9.4/0/-1/{After rotation}}{%
      \begin{scope}[shift={(\ppx,0)}]
        \filldraw[fill=blue!2,draw=blue!14,rounded corners=5pt]
          (-1.65,-1.65) rectangle (1.65,1.70);
        \node[text=ppink,font=\bfseries\small] at (0,2.04) {\pptitle};
        \draw[gray!28,densely dotted] (0,0) circle (1.13);
        \draw[gray!28,densely dotted] (0,0) circle (0.69);
        \fill[gray!55] (0,0) circle (0.025);
        \ifdim\ppx pt>9pt
          \foreach \ppj in {0,...,\numexpr\ppk-1\relax}{%
            \pgfmathsetmacro{\ppinitiala}{360*\ppj/\ppk}
            \pgfmathsetmacro{\ppinitialb}{180/\ppk+360*\ppj/\ppk}
            \draw[green!45!black,line width=0.65pt]
              (\ppinitiala:1.13) circle (0.087);
            \draw[green!45!black,line width=0.65pt]
              (\ppinitialb:0.69) circle (0.087);
          }
        \fi
        \foreach \ppj in {0,...,\numexpr\ppk-1\relax}{%
          \pgfmathsetmacro{\ppa}{\ppphase+360*\ppj/\ppk}
          \pgfmathsetmacro{\ppb}{\ppphase+\ppoffset*180/\ppk+360*\ppj/\ppk}
          \pgfmathsetmacro{\ppanext}{\ppa+360/\ppk}
          \pgfmathsetmacro{\ppbnext}{\ppb+360/\ppk}
          \draw[ppred!40,line width=0.45pt]
            (\ppa:1.13)--(\ppanext:1.13);
          \draw[ppink!35,line width=0.45pt]
            (\ppb:0.69)--(\ppbnext:0.69);
        }
        \foreach \ppj in {0,...,\numexpr\ppk-1\relax}{%
          \pgfmathsetmacro{\ppa}{\ppphase+360*\ppj/\ppk}
          \pgfmathsetmacro{\ppb}{\ppphase+\ppoffset*180/\ppk+360*\ppj/\ppk}
          \pgfmathtruncatemacro{\ppfirstlabel}{\ppj+1}
          \pgfmathtruncatemacro{\ppsecondlabel}{\ppj+\ppk+1}
          \draw[->,ppred,line width=0.9pt] (\ppa:1.13)
            --++({\ppa+90}:0.28);
          \draw[->,ppink,line width=0.9pt] (\ppb:0.69)
            --++({\ppb-90}:0.23);
          \filldraw[fill=ppred,draw=white,line width=0.5pt]
            (\ppa:1.13) circle (0.052);
          \filldraw[fill=ppink,draw=white,line width=0.5pt]
            (\ppb:0.69) circle (0.052);
          \node[text=ppred,font=\bfseries\small] at (\ppa:1.40)
            {\ppfirstlabel};
          \node[text=ppink,font=\bfseries\small] at (\ppb:0.94)
            {\ppsecondlabel};
        }
      \end{scope}
    }
    \draw[->,ppblue,line width=1pt] (1.83,0.18)--(2.87,0.18);
    \node[text=ppblue,font=\footnotesize,align=center] at (2.35,0.65)
      {motion\\to $T$};
    \draw[->,ppblue,line width=1pt] (6.53,0.18)--(7.57,0.18);
    \node[text=ppblue,font=\small] at (7.05,0.68) {$R_{-\theta_1}$};
    \node[text=ppblue,font=\footnotesize,align=center] at (7.05,-0.35)
      {rotate\\back};
    \ifnum\ppk=2
      \node[text=ppink,font=\footnotesize,align=center] at (4.7,-2.08)
        {After rotation: bodies $1,2$ return to their initial states;\\
         body $3$ occupies the initial state of $4$, and $4$ that of $3$.};
    \else
      \node[text=ppink,font=\footnotesize,align=center] at (4.7,-2.08)
        {After rotation: bodies $1,2,3$ return to their initial states;\\
         bodies $4,5,6$ occupy the initial states of $6,4,5$, respectively.};
    \fi
  \end{scope}
}
\end{tikzpicture}
\endgroup}
\caption{The return up to rotation and relabeling for $k=2$ and $k=3$.
Red denotes the first group and black the second; the numbers are body
labels. In each row, rotating the state at $T$ by $-\theta_1$ recovers
the initial state after relabeling the black bodies. Arrows attached to
the bodies indicate tangential velocities, not trajectories. The radii,
angular advance, and arrow lengths are schematic. Green rings in the
last column mark the initial positions; the rotated dots lie exactly
on these rings, with the indicated change of labels.}
\label{fig:pseudoperiodicreturn}
\end{figure}

In our earlier work, we used the equations at $t=T$
in~\eqref{systemneqgeneral} for the numerical search. Here we use
conditions at $H=T/2$ that guarantee~\eqref{systemneqgeneral} at $t=T$.
These are the conditions
\begin{equation}\label{midpointconditions}
\dot r_1(H)=\dot r_2(H)=0,\qquad
\theta(H)-\beta(H)=\frac{\pi}{k},\qquad
\theta(H)=\frac{\theta_1}{2}.
\end{equation}
If $m_2$ and
$\theta_1$ are fixed, these are four equations in the four unknowns
$x_2,x_3,x_4,H$. The first three equations give the reflection symmetry
that we need. The last one fixes the angular advance at $T=2H$.

Notice that the polar angle of body $1$ is $\theta(t)$, whereas
that of body $k+1$ is $\beta(t)+\pi/k$. At $t=H$, the midpoint
conditions give
\[
\theta(H)=\beta(H)+\frac{\pi}{k}=\frac{\theta_1}{2}.
\]
Thus the origin, body $1$, and body $k+1$ are aligned. More precisely,
the two bodies lie on the same ray from the origin, making an angle
$\theta_1/2$ with the positive horizontal axis. Let $\ell_\varphi$ be
the line containing this ray, where $\varphi=\theta_1/2$, and let
$Q_\varphi$ denote reflection in this line. Together with the zero radial velocities at $H$, this alignment allows
us to show that the motion is symmetric with respect to $\ell_\varphi$. 
We prove this rigorously in
Lemma~\ref{midpointlemma}, in Subsection~\ref{midpointsymmetry}, and
show that all eight conditions~\eqref{systemneqgeneral} follow.
These are sufficient conditions for the symmetric motions considered
here.

We use $\theta_1$ as a continuation parameter. For each prescribed
value of $\theta_1$ and a chosen mass ratio $m_2$, we solve the four
midpoint equations~\eqref{midpointconditions} for $x_2,x_3,x_4,H$
using Newton's method. To obtain a solution for a nearby value of
$\theta_1$, we use the previously computed solution as the initial
guess for the next Newton correction.

The stochastic method was used to find initial approximations
close enough for Newton's method to converge. Once such an approximation
is available, Newton's method provides the correction, and continuation
allows us to follow the solutions as $\theta_1$ varies.

The mass ratio may also be varied between computations. In particular,
the examples reported in the tables use several mass ratios; the tables
alone do not establish that all the examples belong to one connected
family, or that the midpoint Jacobian is nonsingular along an entire
branch.

We report selected examples with $\pi/6\leq\theta_1\leq2\pi$ for the
4-body problem and $\pi/6\leq\theta_1\leq\pi$ for the 6-body problem.
The initial conditions in the tables satisfy all eight return equations
\eqref{systemneqgeneral} with numerical residuals smaller than
$10^{-7}$. The tables and verification results are given in
Section~\ref{IC}.

Eight-second videos illustrating the types of motion considered here are available
on the author's YouTube channel. For the 4-body problem, the links
\href{https://youtube.com/shorts/1Y2fvROWLA4}{$30^\circ$},
\href{https://youtube.com/shorts/Cw48LoqNn84}{$45^\circ$},
\href{https://youtube.com/shorts/16Y9WrwK0Nk}{$60^\circ$},
\href{https://youtube.com/shorts/QQAq2uPW6Io}{$90^\circ$},
\href{https://youtube.com/shorts/Im-YlqDY14Q}{$120^\circ$},
\href{https://youtube.com/shorts/_Hi0Qri97Aw}{$135^\circ$},
\href{https://youtube.com/shorts/21cg_kFhxsY}{$150^\circ$},
\href{https://youtube.com/shorts/z9Ckaa0EZjY}{$180^\circ$},
\href{https://youtube.com/shorts/qbqz7sXxTpw}{$210^\circ$},
\href{https://youtube.com/shorts/FK8GesqZ5rI}{$225^\circ$},
\href{https://youtube.com/shorts/KAeUImWuQxA}{$240^\circ$},
\href{https://youtube.com/shorts/pbbxuvHSvXo}{$270^\circ$},
\href{https://youtube.com/shorts/ZjGrJOOnh4Y}{$300^\circ$},
\href{https://youtube.com/shorts/AYUknGq4Vmc}{$315^\circ$},
\href{https://youtube.com/shorts/s-swxm-9HK4}{$330^\circ$}, and
\href{https://youtube.com/shorts/JtAPukLbSFo}{$360^\circ$}
lead to the corresponding videos.

For the 6-body problem, the videos are given by the links
\href{https://youtube.com/shorts/yQZMVa1n6xs}{$30^\circ$},
\href{https://youtube.com/shorts/09HlrbwGmAU}{$45^\circ$},
\href{https://youtube.com/shorts/3ToU4R2GPTI}{$60^\circ$},
\href{https://youtube.com/shorts/3QAtj1A3pFE}{$90^\circ$},
\href{https://youtube.com/shorts/JocBIvAM2r0}{$120^\circ$},
\href{https://youtube.com/shorts/zL5IuIMFq2U}{$135^\circ$},
\href{https://youtube.com/shorts/g8zQE8wSgRk}{$150^\circ$}, and
\href{https://youtube.com/shorts/1V-4WEy9mZs}{$180^\circ$}.
These links also appear in the tables. The videos illustrate the shapes
of the motions; the parameters used in this paper are those given in
the tables.

\begin{figure}[t]
\centering
\includegraphics[width=0.95\linewidth]{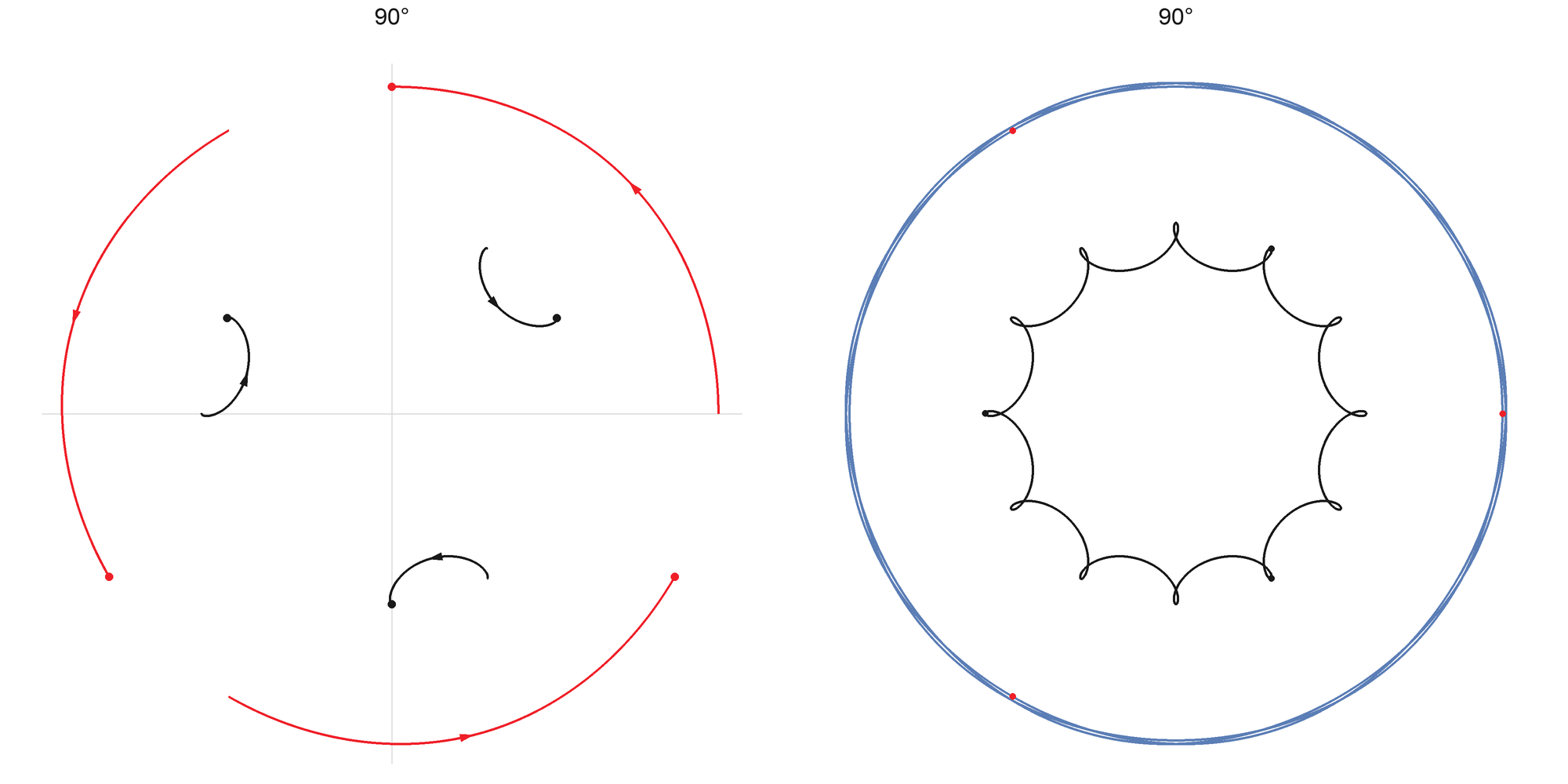}
\caption{The six-body example with $\theta(T)=90^\circ$.
The left image shows all six bodies from $0$ to $T$; the right shows the
union of their trajectories over the labeled period $12T$. Red on the
left and blue on the right denote the first group; black denotes the
second. Dots mark positions at $T$ on the left and at $0$ on the right.}
\label{fig:n6theta90}
\end{figure}

Figure \ref{fig:n6theta90} shows the motion for a
solution with $\theta_1=\frac{\pi}{2}$, and Table \ref{tab:ic-n6} gives
numerical data for the selected examples, all with $r_1(0)=1$.

When
\[
\theta(t)-\beta(t)=c_1
\qquad\text{and}\qquad
\frac{r_1(t)}{r_2(t)}=c_2
\]
for constants $c_1$ and $c_2$, the two groups preserve their relative
orientation and relative size throughout the motion. This corresponds to a
motion generated by a central configuration. In the motions considered
here, the configuration is allowed to change shape with time, and finding
such a motion requires solving differential equations with return
conditions. García-Azpeitia and Ize~\cite{GI} prove global bifurcation results for
periodic motions from polygonal relative equilibria. Their polygonal
model has equal masses on a regular polygon and a central mass. In some
planar branches, the polygon separates into two regular polygons. Here
we impose two polygonal groups directly, without a central body, and
allow their masses, radii, and angular variables to differ.

Earlier variational constructions of symmetric planar motions include
Chen~\cite{Chen2003}, who uses estimates of the action of collision
paths to establish existence. Ferrario and Terracini~\cite{FT2004}
give a general method for obtaining collision-free action minimizers
under symmetry constraints. In particular, their Example~11.7 gives
nonhomographic planar six-body motions consisting of two centered
equilateral triangles, with a similar construction for two regular
$k$-gons when $k$ is odd. The symmetry in this example exchanges the
two groups and requires their masses to be equal. Our six-body examples
have different masses in the two triangles.

Shibayama~\cite{Shibayama2006} studies the equal-mass $2k$-body problem
with two regular $k$-gons at every instant. Building on the earlier
variational results, he proves existence of several distinct solutions
and distinguishes them by studying the configurations up to rotation.
His planar solutions have zero total angular momentum. In our units,
\[
\mathcal L=k\bigl(r_1^2\dot\theta+m_2r_2^2\dot\beta\bigr),
\]
which is nonzero for every example in Tables~\ref{tab:ic-n4}
and~\ref{tab:ic-n6}. Hence these examples are not transformed copies
of his planar solutions under scaling, rotation, relabeling, a time
shift, or time reversal.

Yu~\cite{Yu2018} constructs collision-free spatial double choreographies
with equal masses by minimizing the action under symmetry and topological
constraints. His construction is motivated by rotating regular $k$-gons,
but the resulting bodies in each group are related by time shifts along
a common spatial curve. This does not impose a regular polygon at every
instant. In contrast, our bodies remain in one fixed plane, and each
group is related by simultaneous rotations through multiples of $2\pi/k$,
as in~\eqref{generalansatz}.

For four bodies, Chopovda and Sweatman~\cite{CS2022} compute families
of relative periodic orbits with two opposite pairs of equal masses,
allowing the masses of the two pairs to differ. Their families originate
from collinear Schubart orbits and are computed by continuation and
differential correction. 

For six bodies, Bashir~\cite{Bashir2023} studies two coplanar equilateral
triangles and computes equal-mass families originating from a Schubart
orbit and from the regular hexagonal orbit. Her equations allow different
masses between the triangles, but the computed families have equal
masses; the corresponding unequal-mass families are proposed as further
work. Our six-body examples have different masses in the two triangles.
In the present paper, we use the common polygonal description for both
problems, prescribe the angular advance $\theta_1$, and formulate the
return conditions at $H=T/2$. Lemma~\ref{midpointlemma} justifies these
conditions, which we solve using Newton's method, starting from
approximations obtained by the stochastic search.

The stochastic solver treats the return equations as a black box and uses only function
evaluations. At each iteration, several nearby candidates are sampled and
the search box is updated componentwise using the displacement from the last
successful iteration. The equations themselves are deterministic; the term
stochastic refers only to the random selection of trial points. This
search was used only to obtain nearby initial approximations; Newton's
method is used for the subsequent corrections and continuation. The
stochastic method is related to derivative-free direct-search
methods; see \cite{DRRZ}. Related stochastic methods are discussed
in \cite{LPO}.

Section~\ref{DE} derives the differential equations for $\theta$,
$\beta$, $r_1$, and $r_2$ for general $k\geq2$ and group masses
$1$ and $m_2$. It also gives the explicit four- and six-body cases
and proves the midpoint symmetry lemma. Section \ref{GFS}
describes the adaptive stochastic search used to obtain the initial
approximations. Finally, Section \ref{IC} gives the numerical initial conditions
for the computed examples and illustrates several of the resulting motions.


\section{Reduced equations of motion}\label{DE}
We derive the equations using the radial and tangential orthonormal frames
along the two groups of bodies. Let $k\geq2$ be any integer and
$m_2>0$. Using the ansatz~\eqref{generalansatz}, write
\begin{equation}\label{groupansatz}
\begin{aligned}
{\bf a}_j(t)&=r_1(t)e\bigl(\theta(t)+2\pi j/k\bigr),\\
{\bf b}_j(t)&=r_2(t)e\bigl(\beta(t)+\pi/k+2\pi j/k\bigr),
\qquad j=0,\ldots,k-1.
\end{aligned}
\end{equation}
The bodies ${\bf a}_j$ have mass $1$ and the bodies ${\bf b}_j$ have
mass $m_2$. For $k=2$, these are bodies $1,2$ and $3,4$, respectively.
For $k=3$, they are bodies $1,2,3$ and $4,5,6$.

For convenience, let us define
\begin{equation}\label{distances}
\begin{gathered}
\delta=\theta-\beta,\qquad
\psi_j=\frac{(2j+1)\pi}{k}-\delta,\\
D_j=\sqrt{r_1^2+r_2^2-2r_1r_2\cos\psi_j},\qquad
\lambda_k=\frac14\sum_{j=1}^{k-1}\csc\frac{\pi j}{k}.
\end{gathered}
\end{equation}
Notice that $D_j$ is the distance from ${\bf a}_0$ to ${\bf b}_j$ and that
\[
\lambda_2=\frac14,\qquad\lambda_3=\frac{1}{\sqrt3}.
\]
All computations below take place where $r_1,r_2>0$ and $D_j>0$ for
every $j$. These inequalities exclude collisions within and between the
two groups.

\subsection{The radial and tangential components}
For each $j$, the arguments in~\eqref{groupansatz} are taken as
continuous real-valued angles along the motion. As before, let
$J(x,y)=(-y,x)$. The vectors $e(\theta)$ and $Je(\theta)$ form an
orthonormal basis of the plane. A direct differentiation gives
\begin{equation}\label{polaracceleration}
\ddot{\bf a}_0=
\bigl(\ddot r_1-r_1\dot\theta^2\bigr)e(\theta)
+\bigl(2\dot r_1\dot\theta+r_1\ddot\theta\bigr)Je(\theta).
\end{equation}
Since rotations preserve distances, we can compute the components of the
gravitational acceleration at an instant when $\theta=0$.

The contribution to the acceleration of ${\bf a}_0$ from the other bodies
in the first group is radial. Indeed, the tangential components cancel in
pairs. For a body with angular difference $\gamma=2\pi j/k$, its radial
contribution is
\[
\frac{r_1(\cos\gamma-1)}{(2r_1\sin(\gamma/2))^3}
=-\frac{1}{4r_1^2\sin(\pi j/k)}.
\]
Therefore the total contribution is $-\lambda_k/r_1^2$ in the radial
direction. For $k=2$ this is the attraction of the opposite body. For
$k=3$ it is the sum of the attractions of the other two vertices of the
equilateral triangle.

The radial and tangential components due to the second group are $m_2A_1$
and $m_2B_1$, where
\begin{equation}\label{ABfirst}
A_1=\sum_{j=0}^{k-1}\frac{r_2\cos\psi_j-r_1}{D_j^3},\qquad
B_1=\sum_{j=0}^{k-1}\frac{r_2\sin\psi_j}{D_j^3}.
\end{equation}
For the body ${\bf b}_0$, use the orthonormal basis
$e(\beta+\pi/k),Je(\beta+\pi/k)$. The contribution from its own group
is $-m_2\lambda_k/r_2^2$ in the radial direction. The contribution from
the first group has components
\begin{equation}\label{ABsecond}
A_2=\sum_{j=0}^{k-1}\frac{r_1\cos\psi_j-r_2}{D_j^3},\qquad
B_2=-\sum_{j=0}^{k-1}\frac{r_1\sin\psi_j}{D_j^3}.
\end{equation}
To see the sign in $B_2$, notice that the angles from ${\bf b}_0$ to the
first polygon, measured relative to its radial direction, are
$\delta-\pi/k+2\pi j/k$. As $j$ ranges over the indices, these are
$-\psi_j$ up to a permutation and multiples of $2\pi$.

\begin{proposition}\label{reducedODE}
Assume that $G=1$, that the masses are as above, and that the bodies move
according to \eqref{groupansatz} without collisions. Then they solve the
$2k$-body problem if and only if
\begin{equation}\label{ODE}
\begin{aligned}
\ddot r_1&=r_1\dot\theta^2-\frac{\lambda_k}{r_1^2}+m_2A_1,
&\ddot\theta&=\frac{m_2B_1-2\dot r_1\dot\theta}{r_1},\\
\ddot r_2&=r_2\dot\beta^2-\frac{m_2\lambda_k}{r_2^2}+A_2,
&\ddot\beta&=\frac{B_2-2\dot r_2\dot\beta}{r_2},
\end{aligned}
\end{equation}
where $A_i,B_i$ are given by \eqref{ABfirst} and \eqref{ABsecond}.
\end{proposition}
\begin{proof}
Comparing the two components in \eqref{polaracceleration} with the
components of the gravitational acceleration gives the equations for
$r_1$ and $\theta$. The same computation in the frame along ${\bf b}_0$
gives the other two equations. Each remaining body is obtained from one
of these bodies by a rotation through $2\pi j/k$. Since this rotation
preserves distances and the masses within each group are equal, its
equation follows by applying the same rotation. This also proves the
converse.
\end{proof}

\subsection{The two cases explicitly}
For $k=2$, let
\[
d_1=\sqrt{r_1^2+r_2^2-2r_1r_2\sin\delta},\qquad
d_2=\sqrt{r_1^2+r_2^2+2r_1r_2\sin\delta}.
\]
Then $D_0=d_1,D_1=d_2$, and the four components are
\begin{equation}\label{components4}
\begin{aligned}
A_1&=\frac{r_2\sin\delta-r_1}{d_1^3}
    +\frac{-r_2\sin\delta-r_1}{d_2^3},
&B_1&=r_2\cos\delta\left(\frac1{d_1^3}-\frac1{d_2^3}\right),\\
A_2&=\frac{r_1\sin\delta-r_2}{d_1^3}
    +\frac{-r_1\sin\delta-r_2}{d_2^3},
&B_2&=r_1\cos\delta\left(\frac1{d_2^3}-\frac1{d_1^3}\right).
\end{aligned}
\end{equation}
Substituting these expressions and $\lambda_2=1/4$ in \eqref{ODE} gives
the four differential equations for the 4-body problem.

For $k=3$, define
\[
\begin{aligned}
d_1&=\sqrt{r_1^2+r_2^2-2r_1r_2\sin(\pi/6+\delta)},\\
d_2&=\sqrt{r_1^2+r_2^2+2r_1r_2\cos\delta},\\
d_3&=\sqrt{r_1^2+r_2^2-2r_1r_2\sin(\pi/6-\delta)}.
\end{aligned}
\]
Here $D_0=d_1,D_1=d_2,D_2=d_3$. A direct computation gives
\begin{equation}\label{components6}
\begin{aligned}
A_1&=\frac{r_2\sin(\pi/6+\delta)-r_1}{d_1^3}
     +\frac{-r_2\cos\delta-r_1}{d_2^3}
     +\frac{r_2\sin(\pi/6-\delta)-r_1}{d_3^3},\\
B_1&=\frac{r_2\cos(\pi/6+\delta)}{d_1^3}
     +\frac{r_2\sin\delta}{d_2^3}
     -\frac{r_2\cos(\pi/6-\delta)}{d_3^3},\\
A_2&=\frac{r_1\sin(\pi/6+\delta)-r_2}{d_1^3}
     +\frac{-r_1\cos\delta-r_2}{d_2^3}
     +\frac{r_1\sin(\pi/6-\delta)-r_2}{d_3^3},\\
B_2&=-\frac{r_1\cos(\pi/6+\delta)}{d_1^3}
     -\frac{r_1\sin\delta}{d_2^3}
     +\frac{r_1\cos(\pi/6-\delta)}{d_3^3}.
\end{aligned}
\end{equation}
Substituting these expressions and $\lambda_3=1/\sqrt3$ in \eqref{ODE}
gives the equations for the 6-body problem. Thus the same radial and
tangential computation explains both systems.

\subsection{The midpoint symmetry}\label{midpointsymmetry}

We now prove the assertion about the conditions at $H=T/2$ stated in
the introduction. We use the notation ${\bf a}_j$, ${\bf b}_j$, and
$J(x,y)=(-y,x)$ introduced there.

\begin{lemma}\label{midpointlemma}
Let $k\geq2$ be an integer, and assume that a solution with initial conditions
\eqref{ic} has no collisions on $[0,H]$, where $H>0$. If
\eqref{midpointconditions} holds, then the solution extends without
collisions to $[0,2H]$ and satisfies all eight conditions
\eqref{systemneqgeneral} at $T=2H$.
\end{lemma}

\begin{proof}
Write $\varphi=\theta(H)$. The relative-angle condition gives
$\beta(H)+\pi/k=\varphi$. Thus, at $H$, the vertices in the two
groups lie on the same radial directions:
\[
{\bf a}_j(H)=r_1(H)e\left(\varphi+\frac{2\pi j}{k}\right),\qquad
{\bf b}_j(H)=r_2(H)e\left(\varphi+\frac{2\pi j}{k}\right).
\]
Let $Q_\varphi$ denote reflection in the line through the origin with
direction $e(\varphi)$. In coordinates,
\[
Q_\varphi=
\begin{pmatrix}
\cos(2\varphi)&\sin(2\varphi)\\
\sin(2\varphi)&-\cos(2\varphi)
\end{pmatrix},\qquad
Q_\varphi e(\alpha)=e(2\varphi-\alpha).
\]
Taking the indices modulo $k$, we obtain
\begin{equation}\label{midpointpositionreflection}
Q_\varphi{\bf a}_{-j}(H)={\bf a}_j(H),\qquad
Q_\varphi{\bf b}_{-j}(H)={\bf b}_j(H).
\end{equation}
Notice that this reflection preserves each group of equal masses.

Since $\dot r_1(H)=\dot r_2(H)=0$, the velocities at $H$ are
tangent to their corresponding circles. Consequently,
\[
\dot{\bf a}_j(H)=\dot\theta(H)J{\bf a}_j(H),\qquad
\dot{\bf b}_j(H)=\dot\beta(H)J{\bf b}_j(H).
\]
Reflection reverses a quarter turn, that is,
$Q_\varphi J=-JQ_\varphi$. It follows that
\begin{equation}\label{midpointvelocityreflection}
-Q_\varphi\dot{\bf a}_{-j}(H)=\dot{\bf a}_j(H),\qquad
-Q_\varphi\dot{\bf b}_{-j}(H)=\dot{\bf b}_j(H).
\end{equation}
The minus sign is exactly the change in velocity caused by reversing
time.

For $0\leq u\leq H$, reflect the first half of the motion and traverse
it in reverse order by setting
\[
\widetilde{\bf a}_j(H+u)=Q_\varphi{\bf a}_{-j}(H-u),\qquad
\widetilde{\bf b}_j(H+u)=Q_\varphi{\bf b}_{-j}(H-u).
\]
Reflections preserve distances, and reversing time does not change a
second derivative. Therefore these functions solve the same Newton
equations. By~\eqref{midpointpositionreflection} and
\eqref{midpointvelocityreflection}, their positions and velocities at
$H$ agree with the original ones. Uniqueness for the initial value
problem gives the continuation to $2H$. Its distances between bodies
are those in the first half, up to a relabeling, so no collision occurs.

In terms of the four functions in the ansatz, this continuation reads
\begin{equation}\label{reflectionidentities}
\begin{aligned}
r_i(H+u)&=r_i(H-u),\qquad i=1,2,\\
\theta(H+u)&=2\varphi-\theta(H-u),\\
\beta(H+u)&=2\varphi-\frac{2\pi}{k}-\beta(H-u).
\end{aligned}
\end{equation}
Taking $u=H$ and using~\eqref{ic} and $2\varphi=\theta_1$ gives
\[
\begin{gathered}
r_1(2H)=1,\qquad r_2(2H)=x_2,\\
\theta(2H)=\theta_1,\qquad
\beta(2H)=\theta_1-\frac{2\pi}{k}.
\end{gathered}
\]
Differentiating~\eqref{reflectionidentities}, we also obtain
\begin{equation}\label{midpointvelocityparity}
\begin{gathered}
\dot r_i(H+u)=-\dot r_i(H-u),\qquad i=1,2,\\
\dot\theta(H+u)=\dot\theta(H-u),\qquad
\dot\beta(H+u)=\dot\beta(H-u).
\end{gathered}
\end{equation}
Hence $\dot r_1(2H)=\dot r_2(2H)=0$,
$\dot\theta(2H)=x_3$, and $\dot\beta(2H)=x_4$. These are precisely
the remaining conditions in~\eqref{systemneqgeneral}.
\end{proof}

\section{An adaptive stochastic black-box method for solving equations}\label{GFS}

The stochastic search described below supplies nearby initial
approximations. We first explain the Newton correction used to refine
these approximations.

\subsection{Newton correction of the midpoint equations}\label{Newton}
Fix $k$, $m_2$, and $\theta_1$. For
$z=(x_2,x_3,x_4,H)$, integrate~\eqref{ODE} with initial
conditions~\eqref{ic} and define
\[
\Phi(z)=\left(\dot r_1(H),\dot r_2(H),
\theta(H)-\beta(H)-\frac{\pi}{k},
\theta(H)-\frac{\theta_1}{2}\right).
\]
Thus the midpoint conditions are $\Phi(z)=0$. At a Newton iteration,
we compute a correction $v$ by solving the linear system
\[
D\Phi(z)v=-\Phi(z)
\]
and replace $z$ by $z+\alpha v$. We start with $\alpha=1$ and reduce
it if the integration fails or the residual does not decrease. The
Jacobian can be computed by differentiating the initial value problem
or by numerical differences. In our correction computations we use
centered differences at high working precision. To preserve positivity,
we may use $\log x_2$ and $\log H$ as unknowns in place of $x_2$ and
$H$. A small number of equations does not by itself guarantee an
invertible Jacobian or convergence; the initial approximation must
also be sufficiently close to a suitable root. After correction,
we set $T=2H$ and independently check all eight return equations over
$[0,T]$, as described in Section~\ref{IC}.

\subsection{The stochastic search}
We describe a gradient-free stochastic method for approximately solving a system of equations
\[
h(Z)=\bigl(h_1(Z),\dots,h_q(Z)\bigr)=(0,\dots,0)\in\mathbb{R}^q,
\]
where $Z=(z_1,\dots,z_n)\in\mathbb{R}^n$ is the vector of unknowns. We assume that $h$ is given as a black box: for a proposed input $Z$ we can attempt to evaluate $h(Z)$, but the evaluation may fail (for instance, the procedure may not terminate, or the computation may be undefined). In that case we interpret $Z$ as lying outside the domain of the black-box map.

For any $Z\in\mathbb{R}^n$ we define the error functional
\[
\Err(Z)=
\begin{cases}
\|h(Z)\|_{\infty}=\max_{1\le j\le q}|h_j(Z)|, & \text{if $h(Z)$ is successfully evaluated},\\[2mm]
+\infty, & \text{if $h(Z)$ cannot be evaluated}.
\end{cases}
\]

\subsection*{Parameters}
The method takes as input:
\begin{enumerate}
\item an initial guess $Z_0\in\mathbb{R}^n$;
\item an initial box (search-box) vector $d=(d_1,\dots,d_n)$ with $d_i>0$;
\item a minimum box vector $d_{\min}=(\bar d_1,\dots,\bar d_n)$ with $\bar d_i>0$;
\item a shrinking factor $\rho\in(0,1)$;
\item a scaling factor $c>0$;
\item a target error $e_g>0$;
\item a sample size $N\in\mathbb{N}$ (number of random trials per iteration);
\item a maximum number of iterations $L_{\max}\in\mathbb{N}$.
\end{enumerate}

The procedure is as follows: we maintain a current best point $Z_{\best}$ and its error $e_{\best}=\Err(Z_{\best})$, together with a memory vector $\Delta_{\last}\in\mathbb{R}^n_{\ge 0}$ that records the last accepted displacement (componentwise). Initialize
\[
Z_{\best}=Z_0,\qquad e_{\best}=\Err(Z_0),\qquad \Delta_{\last}=(0,\dots,0).
\]
If $e_{\best}<e_g$, we return this initial point immediately.

For $\ell=1,\dots,L_{\max}$ perform the following steps.

\medskip
\noindent\textbf{(1) Random sampling in an axis-aligned box.}
Generate $N$ independent random candidates
\[
Z^{(j)} = Z_{\best} + \xi^{(j)},\qquad j=1,\dots,N,
\]
where $\xi^{(j)}=(\xi^{(j)}_1,\dots,\xi^{(j)}_n)$ and each coordinate $\xi^{(j)}_i$ is sampled uniformly from $[-d_i,d_i]$. (Equivalently, $Z^{(j)}_i$ is uniform in $(z_{\best,i}-d_i,\;z_{\best,i}+d_i)$.)

Compute the errors
\[
e_j=\Err\!\bigl(Z^{(j)}\bigr)\in[0,\infty],\qquad j=1,\dots,N,
\]
and let $j^\ast$ be an index such that $e_{j^\ast}=\min_{1\le j\le N} e_j$. Define $e_S=e_{j^\ast}$ and $Z^\ast=Z^{(j^\ast)}$.

\medskip
\noindent\textbf{(2) Accept or reject.}

\smallskip
\noindent\emph{Case A (improvement).} If $e_S<e_{\best}$, accept $Z^\ast$:
\[
\Delta_{\last}=\bigl|Z^\ast-Z_{\best}\bigr|\quad\text{(componentwise)},\qquad
Z_{\best}=Z^\ast,\qquad e_{\best}=e_S.
\]
Update the box vector by
\[
d \leftarrow \max\bigl\{\,d,\; c\,\Delta_{\last},\; d_{\min}\,\bigr\},
\]
where the maximum is taken componentwise. If $e_{\best}<e_g$, stop and output $(Z_{\best},e_{\best})$.

\smallskip
\noindent\emph{Case B (no improvement).} If $e_S\ge e_{\best}$, shrink the box but do not let it collapse:
\[
d \leftarrow \max\bigl\{\,\rho\,d,\; c\,\Delta_{\last},\; d_{\min}\,\bigr\},
\]
again componentwise.


\medskip
If the loop reaches $\ell=L_{\max}$ without meeting the target, we output
$(Z_{\mathrm{best}},e_{\mathrm{best}})$ as the best approximation found.

Algorithm~\ref{alg:asbbs} shows pseudocode for the algorithm we just described.

Our solver can be viewed as a stochastic, derivative-free \emph{direct-search} procedure applied to the merit function
\[
\mathrm{Err}(Z)=\|h(Z)\|_{\infty},
\]
with the convention that $\mathrm{Err}(Z)=+\infty$ whenever $Z$ is outside the domain of $h$ (e.g., when the underlying black-box evaluation fails). 
Direct-search methods maintain an incumbent iterate and a step size, sample the objective at finitely many trial points, and accept a new incumbent whenever an improvement is detected; otherwise the incumbent is kept and the step size is reduced \cite{DRRZ}. 
In our setting, the step size is replaced by an \emph{anisotropic box radius} $d=(d_1,\dots,d_n)$, and at each iteration we draw $N$ random candidates inside the box centered at the current best point $Z_{\best}$. 
The distinctive feature of our implementation is the \emph{learning term}
$\Delta_{\last}$, which records the componentwise magnitude of the last
successful displacement.
This term is incorporated through updates of the form
\[
d \leftarrow \max\{\rho\, d,\; c\,\Delta_{\last},\; d_{\min}\}
\qquad\text{(componentwise)},
\]
so the sampling region cannot collapse and its \emph{shape adapts} to the observed successful moves. 
In particular, the term $c\,\Delta_{\last}$ can prevent components of the
search box from shrinking below the scale of the last successful move,
yielding a data-driven box that reflects the local scaling of the problem
while preserving the lower bound $d_{\min}$.
\subsection{Solving the $4$-body problem}
For the stochastic search we use $N=800$, $L_{\max}=300$, $c=0.9$,
and $\rho=0.9$. To write the return map in these variables,
we take $n=5$ and $q=8$. For $Z=(z_1,\ldots,z_5)$, we set $m_2=z_4$
and solve the differential equations in Proposition~\ref{reducedODE}
for $k=2$, with initial conditions
\begin{equation*}
\begin{gathered}
r_1(0)=1,\qquad r_2(0)=z_1,\qquad \theta(0)=\beta(0)=0,\\
\dot r_1(0)=\dot r_2(0)=0,\qquad
\dot\theta(0)=z_2,\qquad \dot\beta(0)=z_3.
\end{gathered}
\end{equation*}
We set $T=z_5$ and define
\begin{equation*}
\begin{aligned}
h(Z)=\big(&r_1(T)-1,\ r_2(T)-z_1,\ \dot r_1(T),\ \dot r_2(T),\\
         &\dot\theta(T)-z_2,\ \dot\beta(T)-z_3,
           \ \theta(T)-\beta(T)-\pi,\ \theta(T)-\theta_1\big).
\end{aligned}
\end{equation*}
Here $z_1,z_4,z_5>0$ and $\theta_1$ is a prescribed multiple of
$\pi/12$ with $\pi/6\leq\theta_1\leq 2\pi$.
Recall that if the differential equation is not defined for the given conditions up to $t=T$, then we consider that $Z$ is not in the domain of $h$.

\subsection{Solving the $6$-body problem}

For the $6$-body problem, we use the same vector
$Z=(z_1,\ldots,z_5)$ and the same initial conditions. We set $m_2=z_4$,
solve the equations in Proposition~\ref{reducedODE} with $k=3$, and take $T=z_5$.
This time the return map is
\begin{equation*}
\begin{aligned}
h(Z)=\big(&r_1(T)-1,\ r_2(T)-z_1,\ \dot r_1(T),\ \dot r_2(T),\\
         &\dot\theta(T)-z_2,\ \dot\beta(T)-z_3,
           \ \theta(T)-\beta(T)-\tfrac{2\pi}{3},\ \theta(T)-\theta_1\big),
\end{aligned}
\end{equation*}
where $\pi/6\leq\theta_1\leq\pi$.

The stochastic search supplies initial approximations. The tables below
give the initial conditions refined using Newton's method, in the units
fixed in the introduction. We independently verify the rounded values printed in
the tables by integrating both the reduced equations and the full
Cartesian equations, as described in Section~\ref{IC}.

\clearpage

\begin{algorithm}[H] \caption{Adaptive stochastic black-box solver } \label{alg:asbbs} 
\begin{algorithmic}[1]
 \Require Initial guess $Z_0\in\mathbb{R}^n$; initial radii $d\in(0,\infty)^n$; minimum radii $d_{\min}\in(0,\infty)^n$; shrink factor $\rho\in(0,1)$; safety factor $c>0$; goal $e_g>0$; samples $N\in\mathbb{N}$; max iterations $L_{\max}\in\mathbb{N}$. 
 \Ensure Approximate solution $Z_{\mathrm{best}}$ and error $e_{\mathrm{best}}$. \State Define the black-box error 
 
 \[ \Err(Z)= \begin{cases} \|h(Z)\|_\infty=\max_{1\le j\le q}|h_j(Z)|, & Z\in\mathrm{Dom}(h),\\ +\infty, & Z\notin\mathrm{Dom}(h). \end{cases} \] 
 
 \State $Z_{\mathrm{best}}\gets Z_0$ 
 \State $e_{\mathrm{best}}\gets \Err(Z_{\mathrm{best}})$ 
 \If{$e_{\mathrm{best}}<e_g$}
 \State \Return $(Z_{\mathrm{best}},e_{\mathrm{best}})$
 \EndIf
 \State $\Delta_{\mathrm{last}}\gets 0\in\mathbb{R}^n$ \State $d\gets \max(d,d_{\min})$ \Comment{componentwise maximum} \For{$\ell=1,2,\dots,L_{\max}$} \For{$j=1,2,\dots,N$} 
 \State Sample $\xi^{(j)}\in\mathbb{R}^n$ with $\xi^{(j)}_i\sim\mathrm{Unif}[-d_i,d_i]$ independently 
 \State $Z^{(j)}\gets Z_{\mathrm{best}}+\xi^{(j)}$ 
 \State $e_j\gets \Err(Z^{(j)})$ \EndFor \State Choose $j^\ast\in\{1,\dots,N\}$ such that $e_{j^\ast}=\min_{1\le j\le N} e_j$
  \State $Z^\ast\gets Z^{(j^\ast)}$, \quad $e_S\gets e_{j^\ast}$ \If{$e_S < e_{\mathrm{best}}$} \Comment{improvement} 
 \State $\Delta_{\mathrm{last}}\gets |Z^\ast-Z_{\mathrm{best}}|$ \Comment{componentwise absolute value} 
 \State $Z_{\mathrm{best}}\gets Z^\ast$ \State $e_{\mathrm{best}}\gets e_S$ \State $d\gets \max(d,\; c\,\Delta_{\mathrm{last}},\; d_{\min})$ \Comment{prevent collapse} \If{$e_{\mathrm{best}}< e_g$}
  \State \Return $(Z_{\mathrm{best}},e_{\mathrm{best}})$ 
  \EndIf \Else \Comment{no improvement}
 \State $d\gets \max(\rho\,d,\; c\,\Delta_{\mathrm{last}},\; d_{\min})$ 
 \EndIf \EndFor 
 \State  \Return $(Z_{\mathrm{best}},e_{\mathrm{best}})$ 
\end{algorithmic} 
\end{algorithm}


\section{Initial conditions and images}\label{IC}

In this section we give 16 four-body examples and 8 six-body examples.
The initial conditions are of the form~\eqref{ic}, with $r_1(0)=1$
for every example in both tables. Clicking an angle in either table
opens an animation on the author's YouTube channel. Each video illustrates
a motion with the indicated angular advance. Its parameter values need
not coincide with those in the table.

The nonconstant parameters are printed to 16 decimal places. We treat
these decimal numbers as exact inputs and integrate the explicit reduced
equations from $0$ to the full time $T$. For the parameters
$(x_2,x_3,x_4,m_2,T)$ in a row, define
\begin{equation}\label{verifiedresidual}
\begin{aligned}
\mathcal R=\big(&r_1(T)-1,\ r_2(T)-x_2,\ \dot r_1(T),\ \dot r_2(T),\\
 &\dot\theta(T)-x_3,\ \dot\beta(T)-x_4,
   \ \theta(T)-\beta(T)-2\pi/k,\ \theta(T)-\theta_1\big).
\end{aligned}
\end{equation}
Let $E_{\rm red}=\max_j|\mathcal R_j|$. We use 50-digit working
precision with 26-digit accuracy and precision goals, and compare the
endpoint with a 45-digit integration with 20-digit goals. We also
integrate the full Cartesian equations, without imposing the polygon
ansatz, using 45-digit working precision and 23-digit goals. Let
$E_{\rm Cart}$ be the largest absolute difference in a Cartesian position
or velocity from the rotated and relabeled initial state described in
the introduction.
We also evaluate the four midpoint conditions independently from the
integration over $[0,T]$ and denote their largest absolute residual by
$E_H$.

Across all 24 rows, the largest computed return residuals satisfy
\[
E_{\rm red}<1.55\times10^{-14},\qquad E_{\rm Cart}<2.95\times10^{-15}.
\]
The midpoint residuals satisfy $E_H<3.11\times10^{-13}$.
The largest difference between the two reduced endpoint computations is
less than $3.18\times10^{-18}$.

These are numerical return residuals and consistency checks. They are
not rigorous bounds on the error of an entire trajectory. The accompanying
data files give the residuals for each row.

\subsection{Initial conditions for the 4-body problem}

\begin{table}[htbp]
\centering
\scriptsize
\setlength{\tabcolsep}{3.5pt}
\renewcommand{\arraystretch}{1.15}
\resizebox{\linewidth}{!}{%
\begin{tabular}{r r r r r r r}
\hline
$\theta(T)$ (deg) & $r_2(0)$ & $\theta'(0)$ & $\beta'(0)$ & $m_2$ & $T$  & $\ell$ \\
\hline
\href{https://youtube.com/shorts/1Y2fvROWLA4}{30} & 0.1674360580061414 & 1.8739488108309801 & -9.0196353179286476 & 1.6103851903290000 & 0.2787446326461474  & 12 \\
\href{https://youtube.com/shorts/Cw48LoqNn84}{45} & 0.2363417213299243 & 1.9124536016375408 & -5.2987177874173311 & 1.6763145957900000 & 0.4082502398744771  & 8 \\
\href{https://youtube.com/shorts/16Y9WrwK0Nk}{60} & 0.3101735618980479 & 2.0851578468841097 & -3.7188791073539404 & 2.0165671564210000 & 0.4965821203675982  & 6 \\
\href{https://youtube.com/shorts/QQAq2uPW6Io}{90} & 0.4896144338645078 & 1.8084629900120226 & -1.2727625998400442 & 1.5658646583060000 & 0.8472253023359529  & 4 \\
\href{https://youtube.com/shorts/Im-YlqDY14Q}{120} & 0.7361372167068833 & 1.9925218035866452 & -0.5868880119054691 & 2.1397249885460000 & 1.0118241025080113  & 6 \\
\href{https://youtube.com/shorts/_Hi0Qri97Aw}{135} & 0.9106047277758632 & 1.8780205157388285 & -0.2690680989486511 & 2.1281190370420000 & 1.2075368987220102  & 8 \\
\href{https://youtube.com/shorts/21cg_kFhxsY}{150} & 1.1403071969071504 & 1.5039714231665109 & 0.0129964867191455 & 1.6590987844550000 & 1.6823057762352595  & 12 \\
\href{https://youtube.com/shorts/z9Ckaa0EZjY}{180} & 1.9007900367591827 & 0.9425224513080898 & 0.1585495817631031 & 1.1307945638480000 & 3.3784766300537913  & 2 \\
\href{https://youtube.com/shorts/qbqz7sXxTpw}{210} & 3.1016472572496311 & 0.7069584390595789 & 0.1250399606153320 & 1.1694059263590000 & 6.0115326079930815  & 12 \\
\href{https://youtube.com/shorts/FK8GesqZ5rI}{225} & 3.5663377602137075 & 0.5993437466670028 & 0.1103513898224311 & 0.8517558687858732 & 7.8064317001451942  & 8 \\
\href{https://youtube.com/shorts/KAeUImWuQxA}{240} & 3.1290652445654781 & 0.5335944592184383 & 0.1467393099303779 & 0.4696453183093701 & 7.8162125790436007  & 6 \\
\href{https://youtube.com/shorts/pbbxuvHSvXo}{270} & 3.1009831006156736 & 0.3744374990159622 & 0.2700830615932694 & 0.9522863887070000 & 5.9572498636267843  & 4 \\
\href{https://youtube.com/shorts/ZjGrJOOnh4Y}{300} & 2.2664556139715240 & 0.3516960701299723 & 0.4095238576067864 & 0.7648737009120000 & 5.2769288252641579  & 6 \\
\href{https://youtube.com/shorts/AYUknGq4Vmc}{315} & 2.2338573128369682 & 0.3661567779011885 & 0.4101125378613508 & 0.3666343531860000 & 6.0477713288221224  & 8 \\
\href{https://youtube.com/shorts/s-swxm-9HK4}{330} & 1.8663774708402239 & 0.3496551028408331 & 0.5113591467623399 & 0.4571565246130000 & 5.2660562765248354  & 12 \\
\href{https://youtube.com/shorts/JtAPukLbSFo}{360} & 1.5696255452683994 & 0.3202342468104805 & 0.5976925265584344 & 0.2712937235040000 & 5.0435929577357010  & 2 \\
\hline
\end{tabular}%
}
\caption{Initial conditions and return time $T$ for the 4-body examples; $\ell T$ is a labeled period. Click an angle to watch an animation on the author's YouTube channel.}
\label{tab:ic-n4}
\end{table}

\subsection{Initial conditions for the 6-body problem}

\begin{table}[htbp]
\centering
\scriptsize
\setlength{\tabcolsep}{3.5pt}
\renewcommand{\arraystretch}{1.15}
\resizebox{\linewidth}{!}{%
\begin{tabular}{r r r r r r r}
\hline
$\theta(T)$ (deg) & $r_2(0)$ & $\theta'(0)$ & $\beta'(0)$ & $m_2$ & $T$  & $\ell$ \\
\hline
\href{https://youtube.com/shorts/yQZMVa1n6xs}{30} & 0.1575682473940597 & 0.8603239758445557 & -2.4322606576615049 & 0.0533588789610000 & 0.6085431401052875  & 12 \\
\href{https://youtube.com/shorts/09HlrbwGmAU}{45} & 0.2953285078963040 & 1.0989969927854905 & -1.5493769286931822 & 0.2000844166650000 & 0.7131962381385923  & 24 \\
\href{https://youtube.com/shorts/3ToU4R2GPTI}{60} & 0.4259902847013917 & 1.3796796885522806 & -1.0129252178600752 & 0.4091869736340000 & 0.7523144520919331  & 6 \\
\href{https://youtube.com/shorts/3QAtj1A3pFE}{90} & 0.5833122915338344 & 1.2464166257015808 & 0.1081264213484286 & 0.3210328209590000 & 1.2355340431122205  & 12 \\
\href{https://youtube.com/shorts/JocBIvAM2r0}{120} & 0.9599608769774455 & 1.2735887900838723 & 0.4783101982972562 & 0.4763235208900000 & 1.5579584710408089  & 3 \\
\href{https://youtube.com/shorts/zL5IuIMFq2U}{135} & 1.3946987827858746 & 1.1189287530034620 & 0.4384091606154191 & 0.5263541813890000 & 1.9882130849913500  & 24 \\
\href{https://youtube.com/shorts/g8zQE8wSgRk}{150} & 1.8073984610766828 & 0.9863941134992799 & 0.3709771058855375 & 0.5317975286220000 & 2.5826824356596353  & 12 \\
\href{https://youtube.com/shorts/1V-4WEy9mZs}{180} & 2.4780240300353072 & 0.8156657256955101 & 0.2619201531438070 & 0.2883496384030000 & 4.0646948925679663  & 6 \\
\hline
\end{tabular}%
}
\caption{Initial conditions and return time $T$ for the 6-body examples; $\ell T$ is a labeled period. Click an angle to watch an animation on the author's YouTube channel.}
\label{tab:ic-n6}
\end{table}

The figures have been generated from the data in the tables.
The return arcs are obtained by integrating~\eqref{ODE}.
For each integer angular advance $d$ in degrees, the rotated and relabeled
return gives a labeled period $\ell T$, where
\[
\ell=\operatorname{lcm}\left(k,\frac{360}{\gcd(360,d)}\right).
\]
Indeed, for a labeled return after $s$ copies of the arc, the first group
requires $s\theta_1\in2\pi\mathbb Z$, while the second requires
$s(\theta_1-2\pi/k)\in2\pi\mathbb Z$. These conditions are equivalent
to the rotation having returned to the identity and $k$ dividing $s$.
The value $\ell$ is their least positive integer solution. This does not
exclude a smaller period that is not an integer multiple of $T$.
To draw the unions over these periods, we rotate the computed return
arcs; the relabeling is applied when following an individual body.
The return residuals above were
checked by independent integrations over $[0,T]$.

The figures obtained by rotating the computed return arcs illustrate
the numerical motions. They are not rigorous proofs of periodicity for
the rounded decimal initial conditions.

\begin{figure}[p]
\centering
\includegraphics[width=0.95\linewidth]{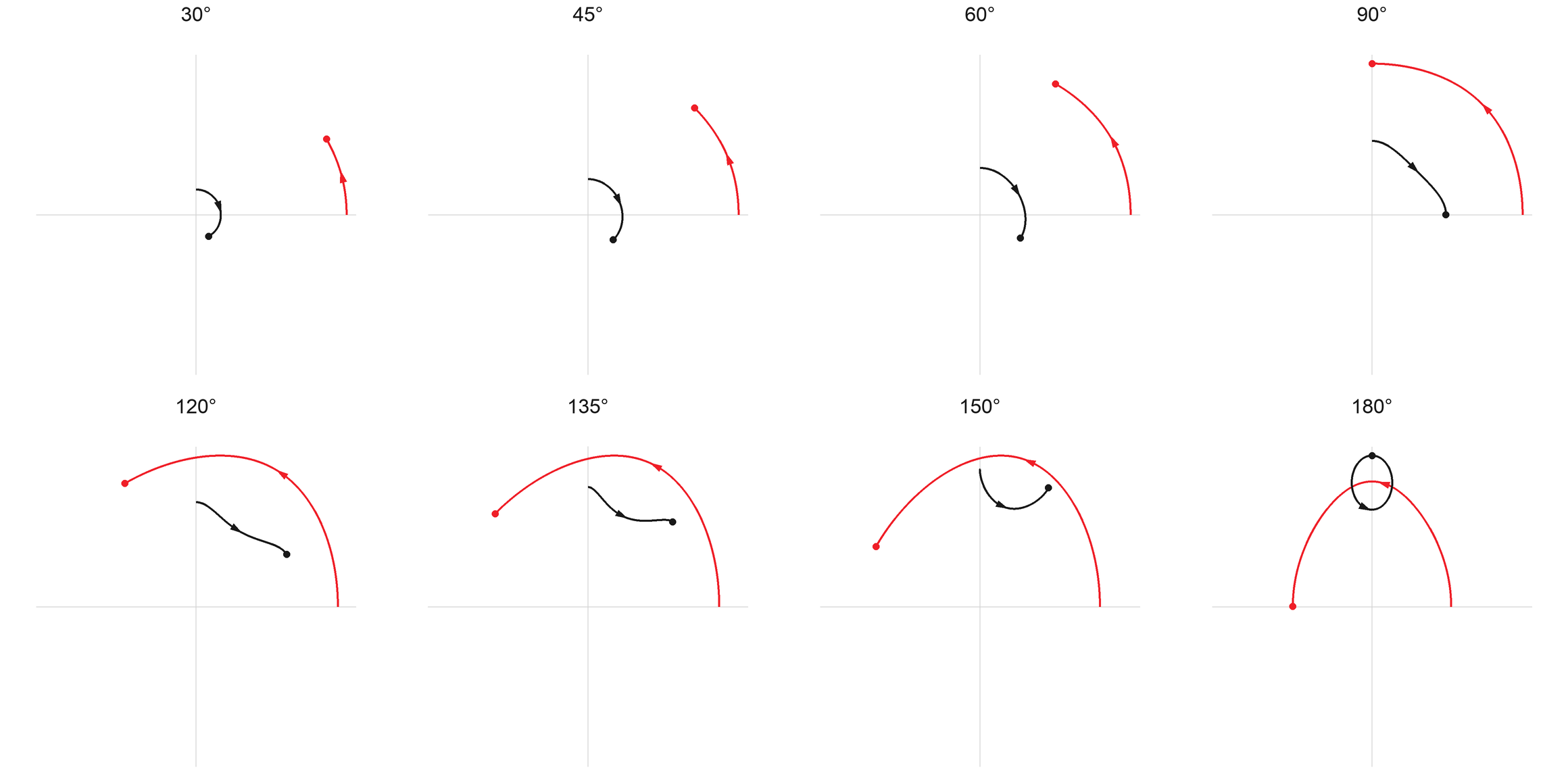}\\[3mm]
\includegraphics[width=0.95\linewidth]{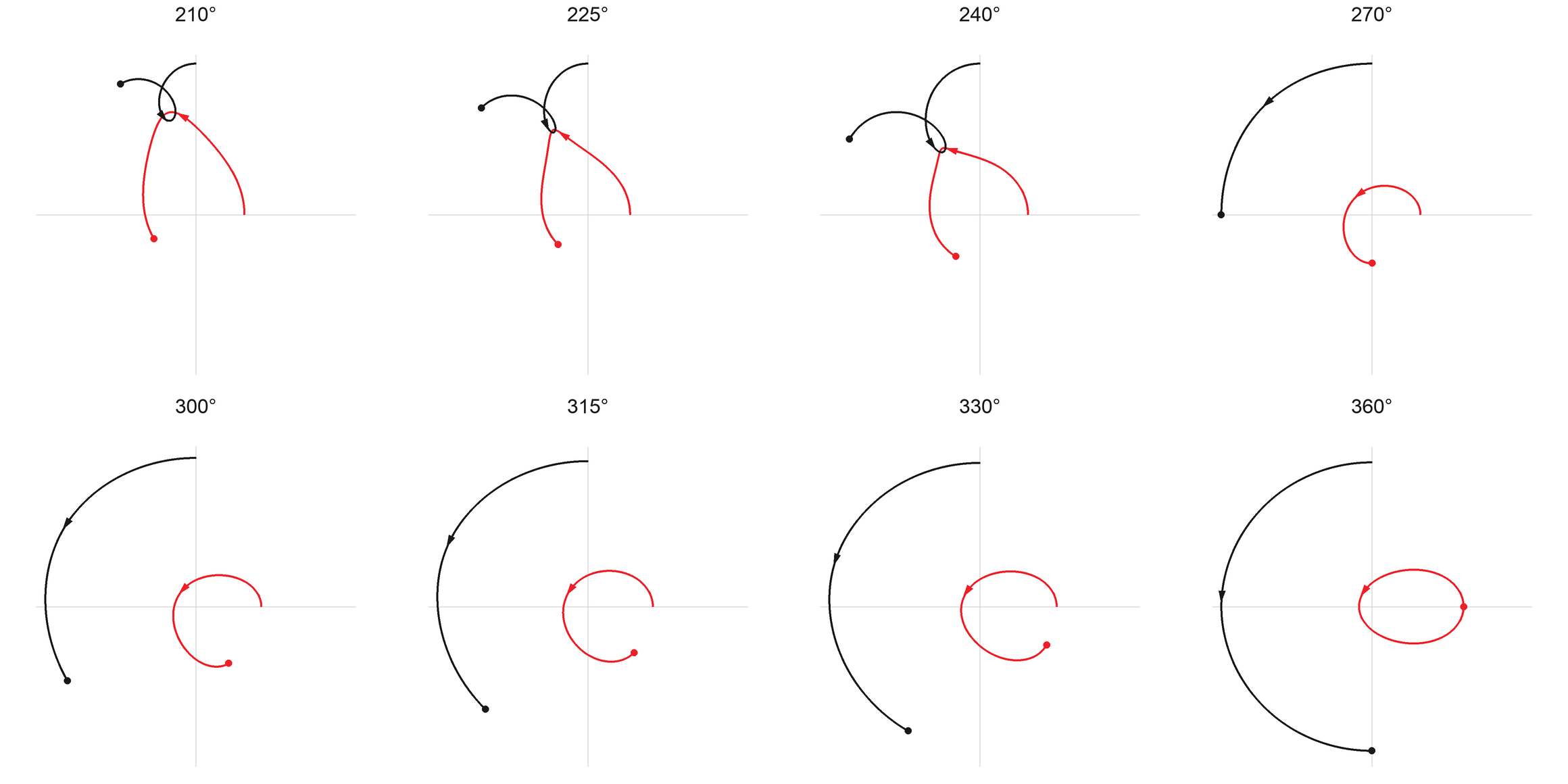}
\caption{Motion of bodies 1 and 3 from $0$ to $T$ for the
four-body examples in Table~\ref{tab:ic-n4}, in the same order. The first
group is red and the second is black. Dots mark positions at $T$ and
arrows show the direction of motion.}
\label{fig:onesteppanel2x8}
\end{figure}

\begin{figure}[p]
\centering
\includegraphics[width=0.95\linewidth]{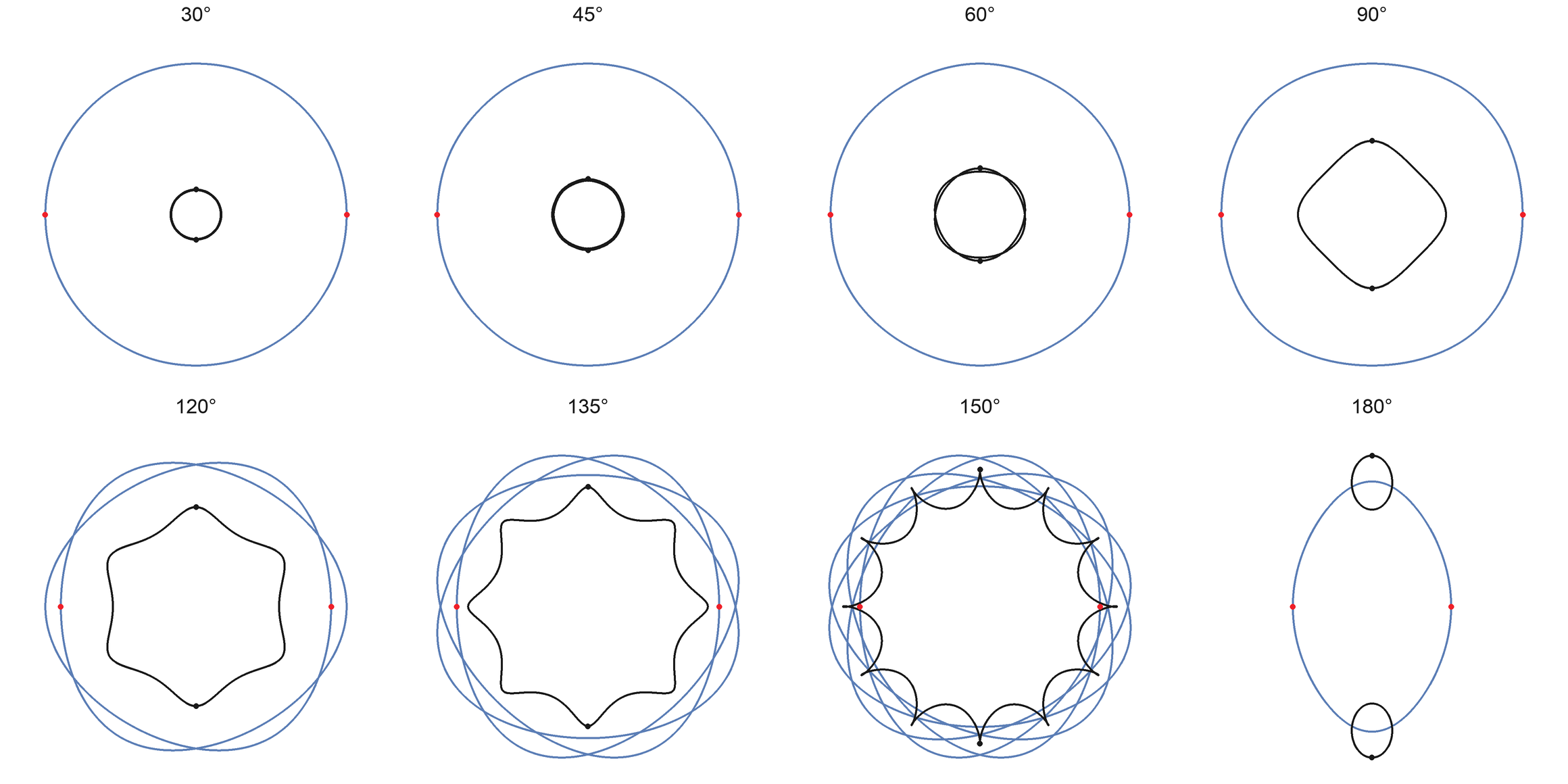}\\[3mm]
\includegraphics[width=0.95\linewidth]{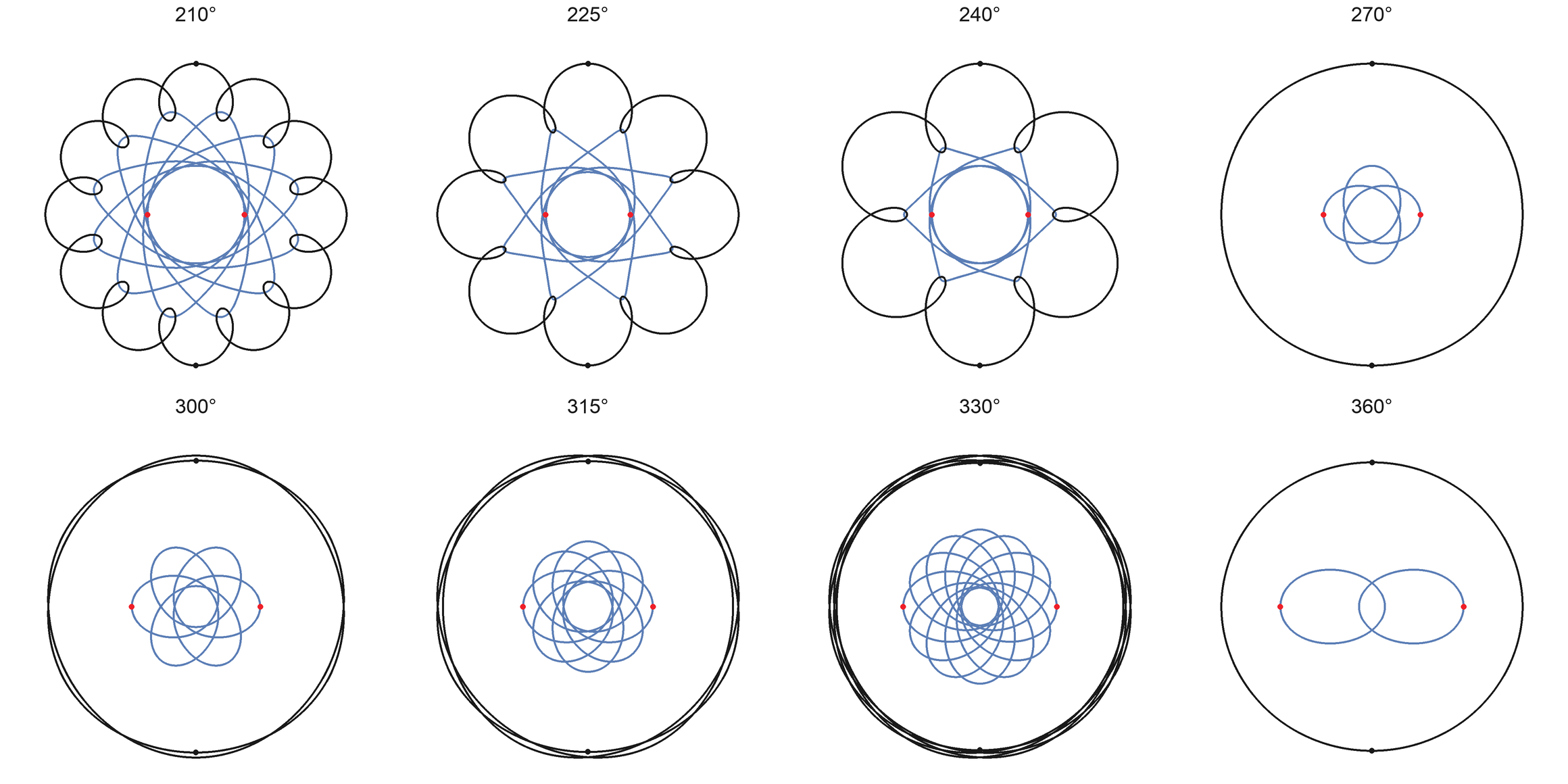}
\caption{Unions of the trajectories over labeled periods for the
four-body examples, in the order of Table~\ref{tab:ic-n4}.
The first group is blue and the second is black. Dots mark positions at
$0$. All panels use equal spatial scales on their two coordinate axes.}
\label{fig:n4panel2x8}
\end{figure}

\begin{figure}[p]
\centering
\includegraphics[width=0.95\linewidth]{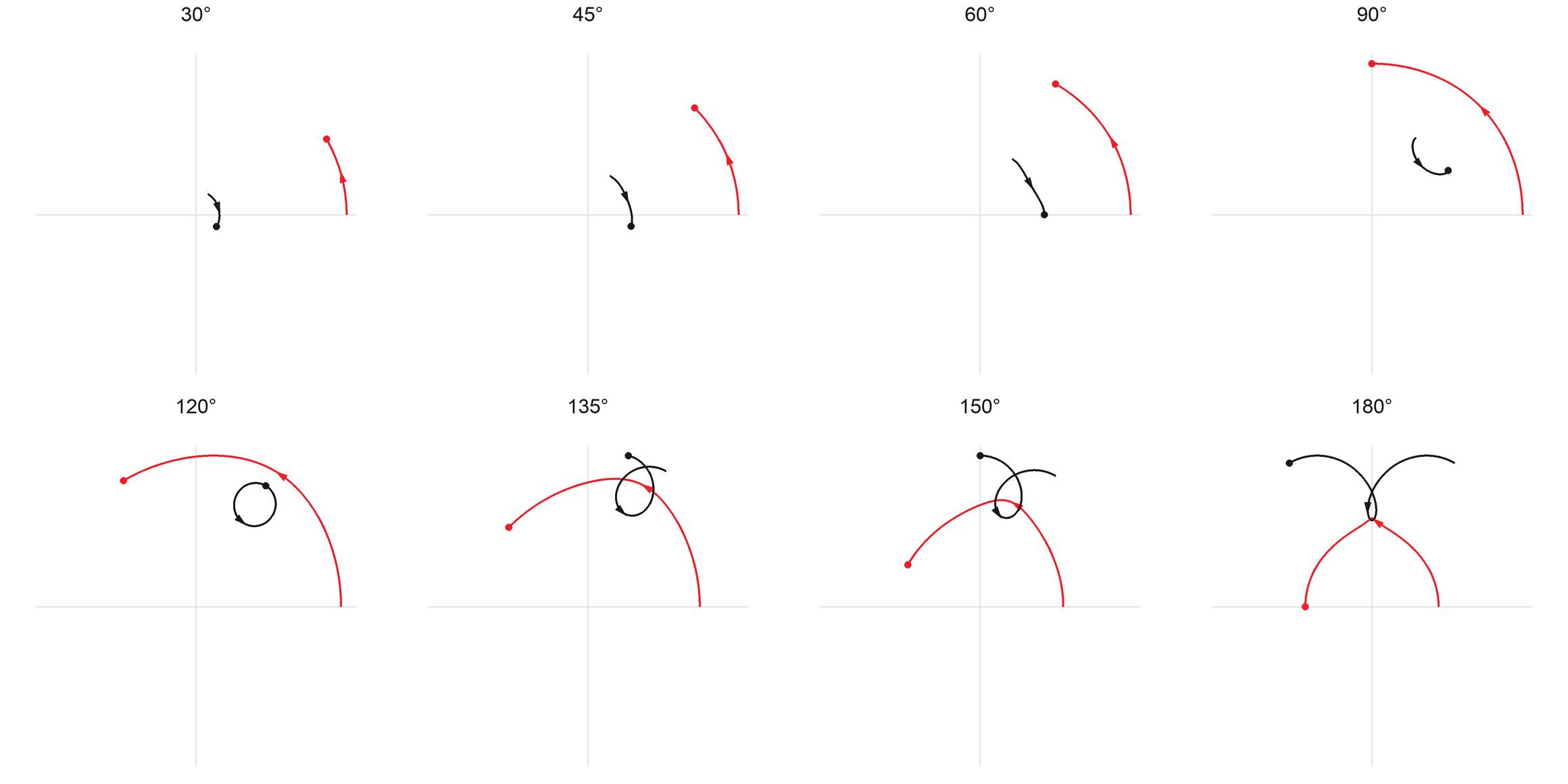}
\caption{Motion of bodies 1 and 4 from $0$ to $T$ for the
six-body examples in Table~\ref{tab:ic-n6}. The first group is red and
the second is black; dots mark positions at $T$.}
\label{fig:onestepn6panel2x8}
\end{figure}

\begin{figure}[p]
\centering
\includegraphics[width=0.95\linewidth]{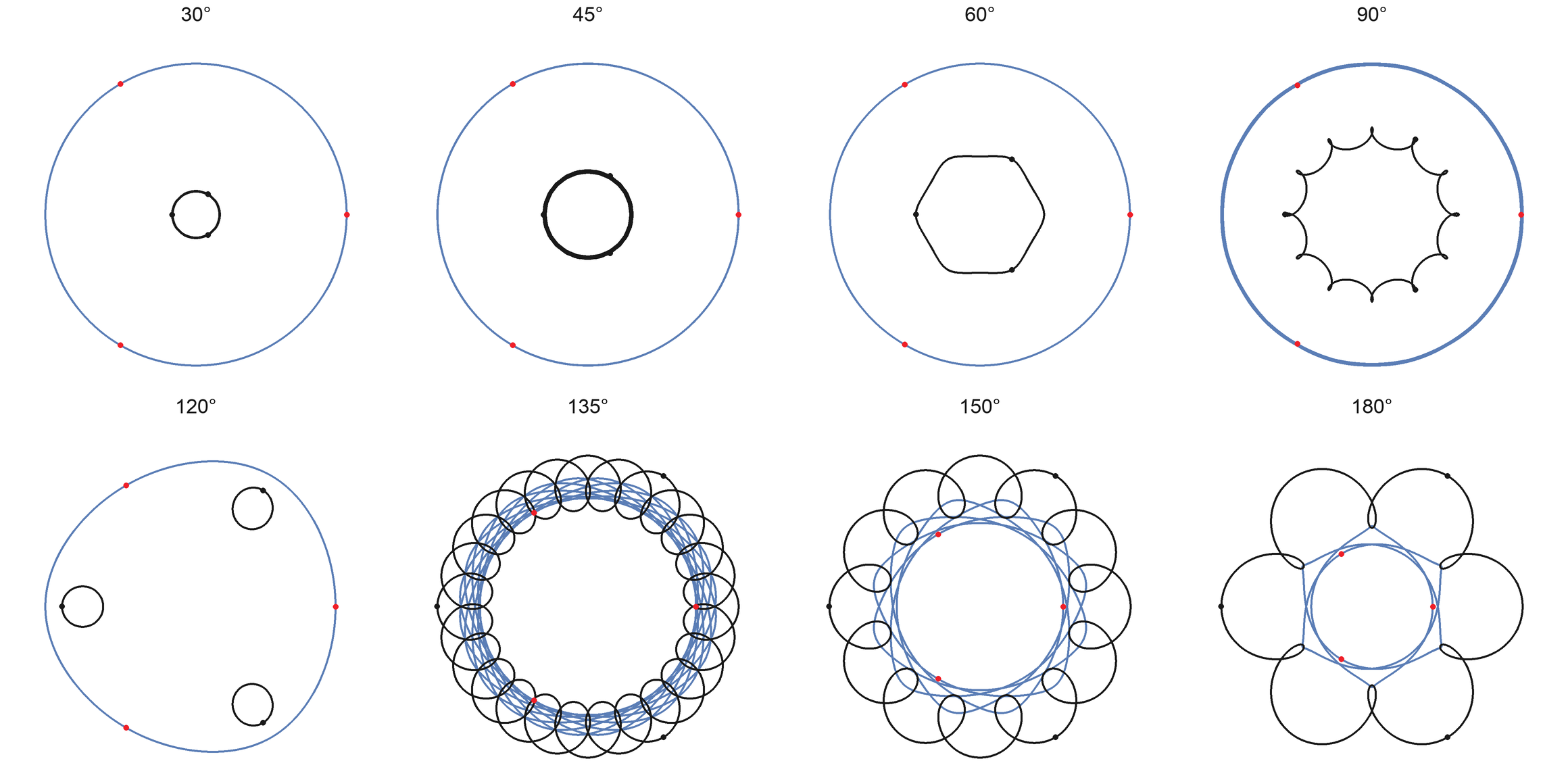}
\caption{Unions of the trajectories over labeled periods for the
six-body examples, in the order of Table~\ref{tab:ic-n6}.
The first group is blue and the second is black; dots mark positions at
$0$.}
\label{fig:n6panel2x8}
\end{figure}
\clearpage

\section*{Declarations}

\textbf{Funding.} This work was supported by an internal AAUP-CSU Research Grant.

\textbf{Competing interests.} The author declares no competing interests.

\textbf{Data availability.} The initial conditions are given in the tables,
in the units fixed in the introduction. The accompanying data files contain
the exact printed inputs, the numerical verification results, and the
sampled trajectories from $0$ to $T$ used for the figures.

\section*{Declaration of generative AI and AI-assisted technologies in the manuscript preparation process}

During the preparation of this work, the author used ChatGPT (OpenAI) to
assist with the readability and organization of the manuscript, the
preparation of numerical code, and mathematical and numerical
checks. The author reviewed and edited the content as needed and takes
full responsibility for the content of the article.

\end{document}